\documentclass[11pt]{article}
\usepackage{amsmath, amsthm, amssymb}    
\usepackage{newtxtext}
 
\usepackage[margin=2.5cm,noheadfoot]{geometry}
\usepackage{graphicx}			
\usepackage[colorlinks=true, urlcolor=blue, citecolor=black, linkcolor=black]{hyperref}  
\usepackage{gensymb}
\usepackage{float}
\usepackage{graphicx,wrapfig,lipsum}
\usepackage{subcaption}			
\usepackage{xfrac}
\usepackage[symbol]{footmisc}
\usepackage{sectsty}
\sectionfont{\centering}

\renewenvironment{thebibliography}[1]{
  \begin{oldthebibliography}{#1}
    \setlength{\itemsep}{0em}
    \setlength{\parskip}{0em}
}
{
  \end{oldthebibliography}
}

\title{Do the Angles of a Triangle Add up to 180$\degree$?
 -- Introducing Non-Euclidean Geometry}

\author{Hanne Kekkonen\textsuperscript{} 
\vspace{10pt}\\
\textsuperscript{}EEMCS, Delft University of Technology, Netherlands; h.n.kekkonen@tudelft.nl} 

\date{}					

\begin{document}

\maketitle

\thispagestyle{empty}

\begin{abstract}
\footnotesize
\noindent How can we convince students, who have mainly learned to follow given mathematical rules, that mathematics can also be fascinating, creative, and beautiful? In this paper I discuss different ways of introducing non-Euclidean geometry to students and the general public using different physical models, including chalksphere, crocheted hyperbolic surfaces, curved folding, and polygon tilings. Spherical geometry offers a simple yet surprising introduction to the topic, whereas hyperbolic geometry is an entirely new and exciting concept to most. Non-Euclidean geometry demonstrates how crafts and art can be used to make complex mathematical concepts more accessible, and how mathematics itself can be beautiful, not just useful. 
\end{abstract}

\section*{Do the Angles of a Triangle Add up to 180$\degree$?}

One of the fundamental truths every child learns in school is that the angles of a triangle add up to 180$\degree$. By the time you graduate from high school, you have probably seen at least three proofs (likely without mathematical rigour) of this fact, starting from cutting or folding a paper triangle to show that when the angles meet in one point, the result is a straight angle. But is it really true for every triangle you could possibly draw? Posing this question to students forces them to think outside the box in which they have been encouraged to stay throughout their whole mathematics education. To help them, one can elaborate: is it really a property of the triangle or perhaps a property of the sheet of paper or chalkboard you drew the triangle on? 

\begin{wrapfigure}{l}{5.7cm}
\includegraphics[width=5.7cm]{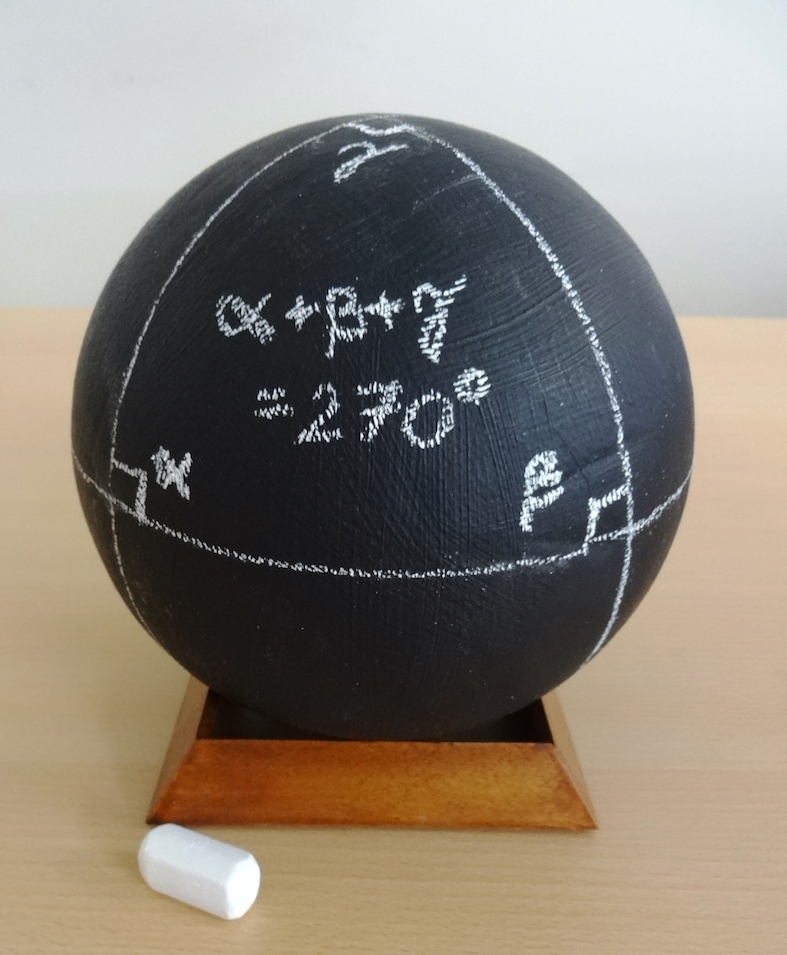}
\caption{A chalksphere.}
\label{Fig:Spherical} 
\end{wrapfigure}
Even though we usually draw and write on flat surfaces, we live on a globe. So, imagine that Alice and Bob are stationed at the north pole and Alice travels directly to the equator. Next, Bob starts his journey at a 90$\degree$ angle to Alice, travels directly to the equator and then follows the equator to Alice. All lines of longitude meet the equator at a 90$\degree$ angle, which means that the created triangle has three right angles. Since we are familiar with the idea of a globe, this example is easy to understand, even for a younger audience, but still offers a surprise for most. You can then ponder with the students; what do straight lines on spheres look like? The simplest way of seeing that straight lines are described by the great circles, that divide the sphere in two equal parts, is to consider a straight line segment as the shortest way between two points. Some other differences between the Euclidean geometry, which we learn in school, and the spherical geometry are that there are no parallel lines, and that the circumference of any circle is 
less than $2\pi r$. A chalksphere or globe offers a great tool for the students to explore spherical geometry, see Figure \ref{Fig:Spherical}\footnote[1]{All the models and pictures presented in this article are by the author.}.

\clearpage
The next step is to guide the students to think about the difference between a sphere and a sheet of paper, and help them arrive at the conclusion that, while paper is flat, a sphere curves. In mathematics this is called curvature. If you draw a tilde $\sim$  you notice that one part of the curve looks like a hill while another part looks like a valley. To distinguish hills from valleys we can assign one of them with positive and the other one with negative curvature. The above example describes one dimensional curvature but how does this translate to two dimensions? If you place your finger on a ball you notice that the ball curves away from it in every direction. If you did the same inside the ball you would notice that there the ball curves towards your hand in every direction. This is why we say that spheres have positive curvature. Comparing an American football and a perfectly round ball we notice that the first one is more curved at the tips than middle, while the second one looks the same everywhere. This means that spheres have constant positive curvature. 

\section*{Hyperbolic Surfaces}
We have now discussed surfaces with zero and positive curvature, which might lead one to wonder: does a surface with negative
curvature exist and how would it look? To attain negative curvature  
the surface would need to curve as a hill in one direction and as a valley in the other, creating a saddle shape 
\begin{wrapfigure}{r}{6.6cm}
\includegraphics[width=6.6cm]{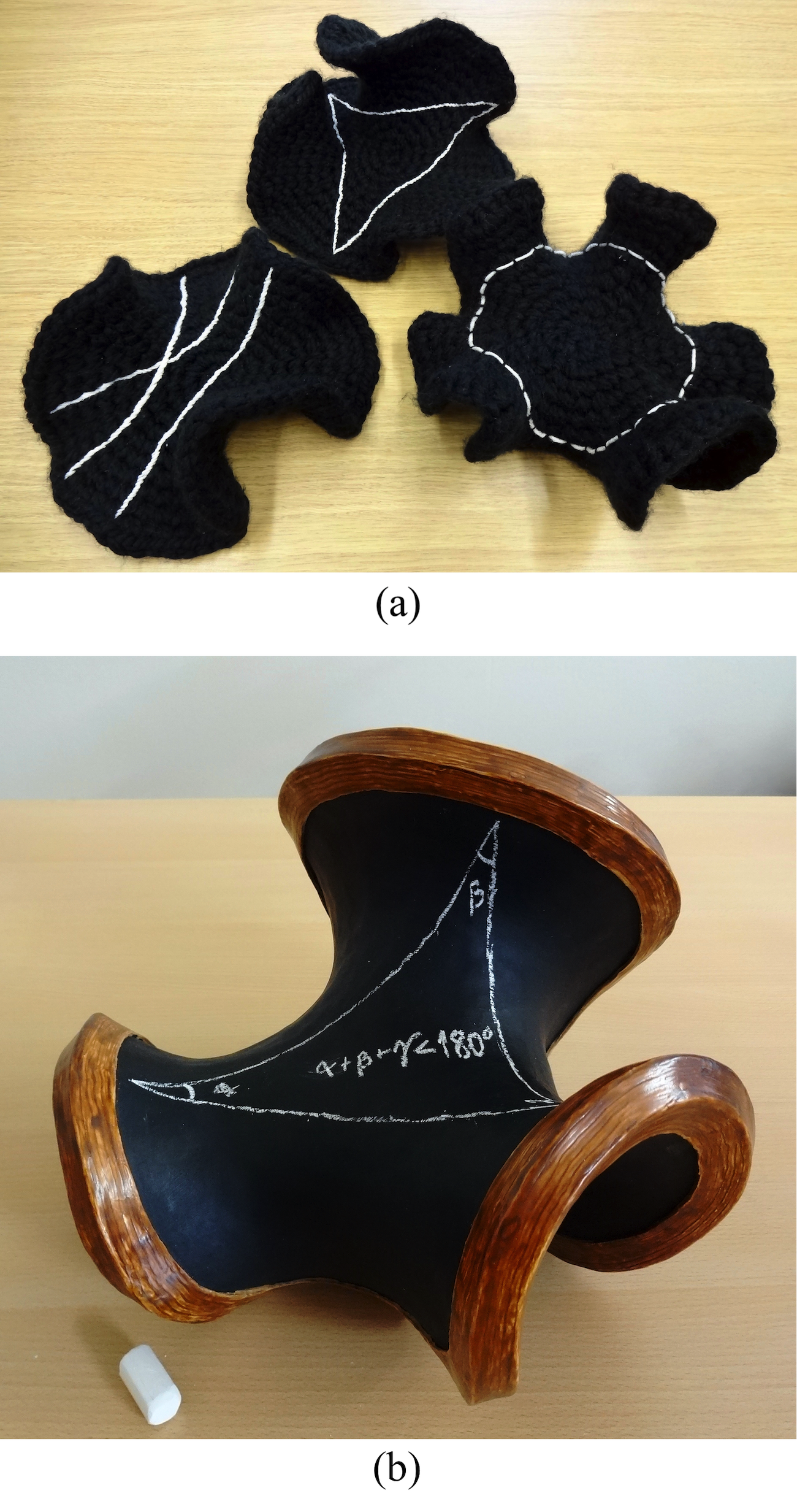}
\vspace{-8mm}
\caption{Three crocheted hyperbolic surfaces (a) and a hyperbolic chalk surface (b).}
\label{Fig:Crochet}
\end{wrapfigure}  
(or a Pringles potato chip). If the curvature is constant we call such a surface hyperbolic. Another way of demonstrating the difference between spheres, flat surfaces and hyperbolic planes is their size. If you try to flatten a half sphere you end up either ripping or stretching it since it covers less area than the plane. Conversely, hyperbolic surfaces are larger than the plane and to flatten one you have to fold it.

Introducing spherical geometry is simple since we are familiar with spheres, but how to introduce hyperbolic surfaces to students? The earliest physical models of hyperbolic surfaces were made by taping identical paper annuli together, see e.g. \cite{Taimina2}. These models are rather fragile and hence would have fairly limited lifespan in the hands of excited members of the audience. Inspired by these paper models, Daina Taimina came up with the idea of crocheting hyperbolic planes to help her differential geometry students \cite{Taimina1,Taimina2}. These models are durable and allow students to explore hyperbolic geometry by playing and folding the surfaces. In Figure \ref{Fig:Crochet} (a) three differences between hyperbolic geometry and Euclidean geometry are demonstrated: the angles of a triangle add up to less than 180$\degree$, a line has infinitely many parallel lines through a given point and the circumference of a circle is more than $2\pi r$. 
I used a crocheted hyperbolic surface as a starting point for creating a hyperbolic chalk surface, shown in Figure \ref{Fig:Crochet} (b). Crocheted models are not only useful and durable, they are also beautiful. 
In addition to Taimina's work, a large art installation curated by Margaret and Christine Wertheim called Hyperbolic Crochet Coral Reef, has been exhibited all around the world.

A second model, that can be used to explain negative curvature, 
is created using a technique called curved folding. Here we start by cutting an annulus out of paper and drawing concentric circles on it. These circles  are then carefully folded with an alternating pattern of mountains and valleys. Paper is an elastic
\begin{figure}[H]
\centering
\begin{minipage}[b]{0.4\textwidth} 
	\includegraphics[height=6.65cm]{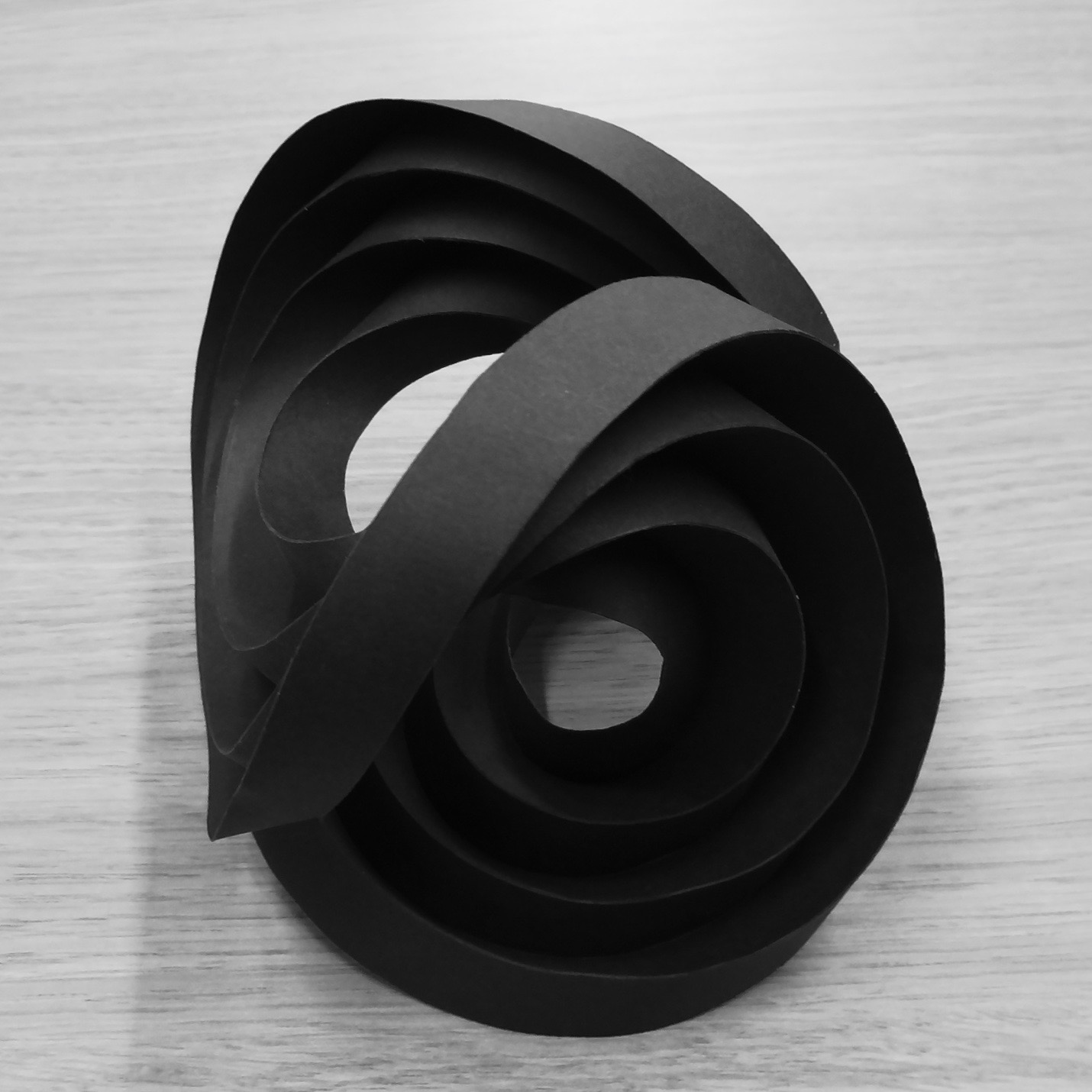}
\end{minipage}
\hspace{10mm}
\begin{minipage}[b]{0.45\textwidth} 
	\includegraphics[height=6.65cm]{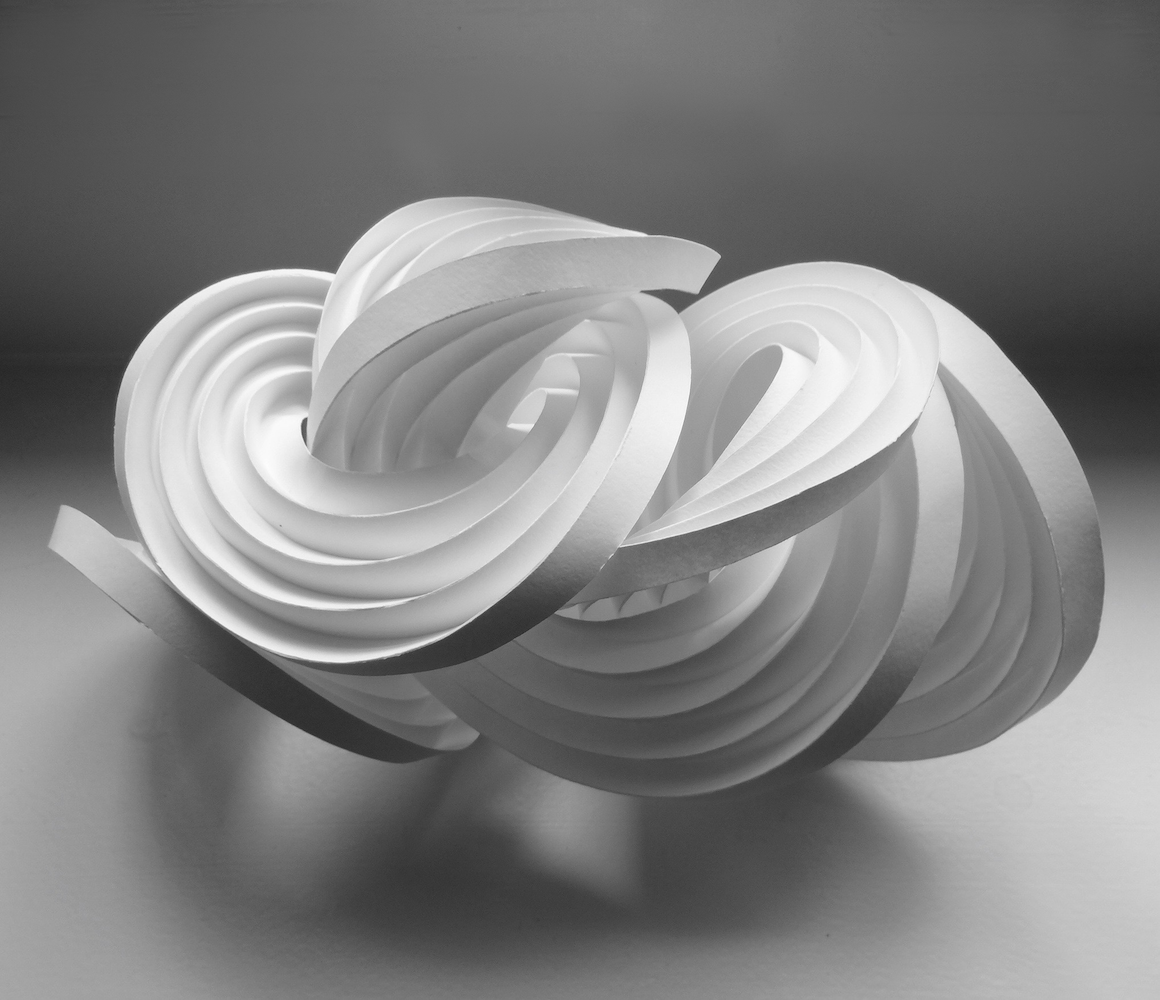}
\end{minipage}
\caption{When a flat paper ring is folded along concentric circles the paper shapes itself into a natural equilibrium form (left). By combining several curved rings one can create paper sculptures.}
\label{Fig:Folding}
\end{figure}
\noindent  material and tries to stay flat unless creased, in which case a plastic deformation occurs. To balance the two forces (staying flat between folds while folded at the creases) the above model automatically folds into a curved 3d-shape as shown in Figure \ref{Fig:Folding} on left. Notice how the circumference of a creased circle stays the same while it moves closer to the origin, making the radius smaller. The mathematical behaviour (or even existence) of these models is still an open question \cite{Demaine, Dias}. The earliest references to curved folding sculptures are from 1927 with a student's work at the Bauhaus. Since then, the technique has fascinated several artists and mathematicians including Thoki Yenn, Kunihiko Kasahara, David Huffman, and Erik and Martin Demaine. Even though the curvature is not necessarily constant, these models offer a different way of understanding negative curvature and allow students to get creative. 

\section*{ Introducing Curvature with Tiling Models}
Another way that can be used to introduce curvature is based on approximating surfaces by tessellations. The honeycomb pattern found in nature consists of hexagons laid next to each other creating a flat surface. Notice that the internal angles of a hexagon are $120\degree$, which means that when three of them meet in a vertex they add up to $360\degree$, creating a flat surface. 
Now, if you start building a surface using a pentagon instead of a  hexagon, \linebreak
\vspace{-3mm}
\begin{figure}[H]
\centering
\begin{minipage}[b]{0.88\textwidth}   
\includegraphics[width=\textwidth]{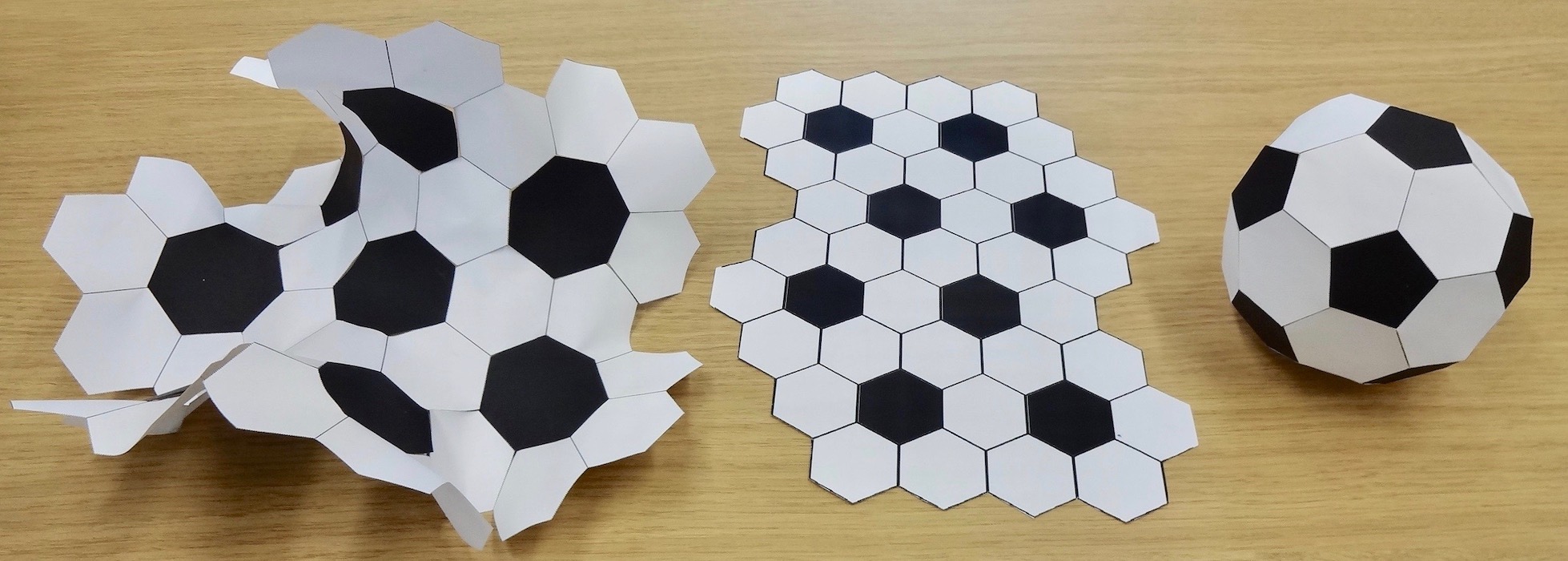}
\end{minipage}
\caption{Polygon tiling models offer a great way of exploring curvature.}
\label{Fig:Tiling}
\end{figure}
\noindent  
and surround it with five hexagons the surface is too small to lie flat and starts to curve. The internal angles of a pentagon are $108\degree$, so the angles at a vertex with two hexagons and one pentagon add up to $348\degree<360\degree$. If continued, this pattern yields a traditional soccer ball approximation of a sphere, shown in Figure \ref{Fig:Tiling}. To create an approximation of a hyperbolic plane, which is too large to lie flat, we start with a heptagon (internal angles $128\sfrac{4}{7}\degree$) and surround it with seven hexagons (the angles at a vertex add up to  $368\sfrac{4}{7}\degree>360\degree$). This method allows people to quickly build a paper model of the hyperbolic plane. It also has the advantage that, while the model is curved, each polygon is made of flat paper and hence it is possible to draw straight lines on it. The model is excellent for more hands-on events allowing students to explore hyperbolic geometry by drawing lines on the models, and seeing how parallel lines diverge and the angles of a triangle add up to less than 180$\degree$. Frank Sottile describes more activities with the ``hyperbolic football'' model explained above \cite{Sottile}.


\section*{ Summary and Conclusions}

Extracurricular mathematics activities usually consist of mathematical puzzles or problems that can be solved by creative use of logic and high school mathematics. Non-Euclidean geometry, on the other hand, introduces an entirely new concept to the audience. It can be introduced without any formulas and gives the students a peek into university level mathematics, which can be very far from the image most people have of mathematics. Spherical geometry is easy to understand, yet offers a couple of surprises, since most people have only encountered Euclidean geometry in school. After introducing the concept of curvature, and the fact that spheres have positive curvature while flat surfaces have zero curvature, one might start to wonder how a surface with negative curvature would look like.
Non-Euclidean geometry offers an excellent opportunity for combining mathematics with crafts and art, and hence breaking some 	preconceptions of both of them. Students can make their own chalkspheres and study how map projections will always distort some areas and force straight lines to look curved. The need for physical models becomes even more obvious when moving to surfaces with negative curvature. The crochet models offer a tactile way of exploring hyperbolic geometry, while curved folding provides a great method for creating paper sculptures. Polygon tiling models allow for studying hyperbolic geometry by drawing on the created hyperbolic plane and more ambitious students can experiment with different hyperbolic tilings to see how they affect the outcome. 

I gave several outreach talks and hands-on workshops in 2019-2020 at the University of Cambridge (UK) introducing non-Euclidean geometry to high school teachers and students of different ages, and organised a stall at the Maths Public Open Day at the Cambridge Science Festival. 
Most of the audience in these events were astonished about how they had never come to think of the fact that triangles on spheres break the ``$180\degree$ rule'' (and several of them seemed rather eager to inform their geometry teacher of this newly found information). They were also fascinated by the crochet and paper models of the hyperbolic plane, and hyperbolic geometry offered a couple of surprises even to many mathematicians. Seeing different models of the hyperbolic plane helped people to realise how the originally rather abstract seeming concept of negative curvature can be observed all around us from Pringles to kale leaves and sea slugs.


{\setlength{\baselineskip}{12pt} 
\raggedright				

} 

\end{document}